\documentclass[11pt]{article}
\usepackage[ansinew]{inputenc}
\usepackage{array}
\usepackage{color}
\usepackage{amsmath,amsthm,amscd,scalerel}
\usepackage{amsxtra}
\usepackage{amstext}
\usepackage{amssymb}
\usepackage{hyperref}
\usepackage{latexsym}
\usepackage{amsfonts,cite}
\usepackage{graphics,epstopdf}
\usepackage{epsfig}
\usepackage{amsfonts}
\usepackage{graphicx}
\usepackage{tcolorbox}
\usepackage{array,tabularx,colortbl}
\usepackage[framemethod=TikZ]{mdframed}

\providecommand{\U}[1]{\protect\rule{.1in}{.1in}}

\begin{document}
	
	\sloppy
	\newtheorem{thm}{Theorem}
	\newtheorem{cor}{Corollary}
	\newtheorem{lem}{Lemma}
	\newtheorem{prop}{Proposition}
	\newtheorem{eg}{Example}
	\newtheorem{defn}{Definition}
	\newtheorem{rem}{Remark}
	\newtheorem{note}{Note}
	\numberwithin{equation}{section}
	
	\thispagestyle{empty}
	\parindent=0mm
	
	\begin{center}
		{\large \textbf{Operational-umbral approach to bivariate degenerate\\ Hermite polynomials and their partial orthogonality}}\\ 					
	
		\vspace{0.20cm}
		
		{\bf Nusrat Raza$^{1}$, Ujair Ahmad$^{2}$, Subuhi Khan$^{3}$}\\
		\vspace{0.15cm}
		{$^{1}$Mathematics Section, Women’s College, Aligarh Muslim University, Aligarh, India.}
		{$^{2,3}$Department of Mathematics, Aligarh Muslim University, Aligarh, India.}\\
		\footnote{$^{*}$This work has been done under Junior Research Fellowship (Ref No. 231610072319 dated:29/07/2023) awarded to the second author by University Grand Commission, New Delhi.}
		\footnote{$^{1}$E-mail:~nraza.math@gmail.com (Nusrat Raza)}
		\footnote{$^{2}$E-mail:~ujairamu1998@gmial.com (Ujair Ahmad)}
		\footnote{$^{3}$E-mail:~subuhi2006@gmail.com (Subuhi Khan)(Corresponding author)}

	\end{center}
	
	\begin{abstract}
		\noindent
The operational calculus associated with special polynomials has proven to be a powerful tool for analyzing and simplifying their properties. This article examines the bivariate degenerate Hermite polynomials with a focus on their differential equations, monomiality properties, operational identities, and partial orthogonality conditions. These polynomials are redeveloped within the framework of umbral formalism, which is further extended to derive certain new results. The study concludes with key observations and insights that highlight the significance of this approach in advancing the understanding of degenerate Hermite polynomials of negative order.
	\end{abstract}
	\parindent=0mm
	\vspace{.25cm}
	
	
	\vspace{0.25cm}
	
	\section{Preliminaries}
Special functions and polynomials are indispensable across various scientific and mathematical domains due to their broad applicability and computational effectiveness. They play a vital role in solving differential equations, modeling physical processes, and providing precise results in mathematics, physics, and engineering. Recently, a number of mathematicians have studied degenerate versions of various special polynomials and corresponding numbers to enhance their analytical properties and broaden their potential applications, see, for example, \cite{KJLK,HR}. A proper understanding of degenerate polynomials requires a well defined series expansion and the corresponding generating function. These polynomials extend classical ones by introducing an additional parameter, say $\lambda$, such that as $\lambda \to 0$, they reduce to their classical counterparts.\\

Degenerate polynomials of two variables offer new tools for analyzing and solving partial differential equations often encountered in physical problems. Based on the definition and generating function of bivariate Hermite polynomials $H_{n}(x,y)$ \cite{APP}, the bivariate degenerate Hermite polynomials $(BVDHP)$ $H_{n}(x,y|\lambda)$ are introduced through the following generating function and series representation \cite{HR}:
\begin{equation}\label{dheq3}
	(1+\lambda)^{(xt+yt^{2})/\lambda} = \sum_{n=0}^{\infty}H_{n}(x,y|\lambda)\dfrac{t^{n}}{n!},
\end{equation}
\begin{equation}\label{dheq14}
	H_{n}(x,y|\lambda) = \sum_{k=0}^{\lfloor n/2 \rfloor}\frac{n!}{(n-2k)!k!}\left[\frac{\log(1+\lambda)}{\lambda}\right]^{n-k}x^{n-2k}y^{k},
\end{equation}
respectively.\\
These polynomials satisfy the generalized heat equation:
\begin{equation}\label{dheq15}
	\frac{\partial}{\partial{y}}H_{n}(x,y|\lambda) = \frac{\lambda}{\log(1+\lambda)}\frac{\partial^{2}}{\partial{x}^{2}}H_{n}(x,y|\lambda)
\end{equation}
with initial condition
\begin{equation}\label{dheq26}
	H_{n}(x,0|\lambda) = \left(\frac{\log(1+\lambda)}{\lambda}\right)^{n}x^{n}
\end{equation}
in the domain $|x|<\infty$, $|y|<\infty$.\\
As $\lambda \to 0$, equation \eqref{dheq3} reduces to the generating function of $H_{n}(x,y)$ \cite{APP}.\\

Solving equation \eqref{dheq15} with initial condition \eqref{dheq26}, the following operational representation for the $BVDHP$ $H_{n}(x,y|\lambda)$ is obtained:
\begin{equation}\label{dheq23}
	H_{n}(x,y|\lambda) = \exp{\left(\dfrac{\lambda}{\log(1+\lambda)}y\partial_{x}^{2}\right)}\left(\dfrac{x\log(1+\lambda)}{\lambda}\right)^{n}.
\end{equation}
Replacement of $x$ by $2x$ and $y$ by $-1$ in generating function \eqref{dheq3} gives
\begin{equation}\label{dheq8}
			(1+\lambda)^{(2xt-t^{2})/\lambda} = \sum_{n=0}^{\infty} H_{n}(2x, -1|\lambda) \frac{t^{n}}{n!} = \sum_{n=0}^{\infty} H_{n}(x|\lambda) \frac{t^{n}}{n!},
\end{equation}
which is the generating function of degenerate Hermite polynomials $H_{n}(x|\lambda)$ \cite{KJLK}.\\
Equation \eqref{dheq8} can be reformulated as:
\begin{equation*}
	\exp\left((2xt - t^{2}) \frac{\log_{e}(1+\lambda)}{\lambda}\right) = \sum_{n=0}^{\infty} H_{n}(x|\lambda) \frac{t^{n}}{n!}.
\end{equation*}
Simplifying above equation, the following relation is obtained:
\begin{equation}\label{dheq59}
	H_{n}\left(x\sqrt{\frac{\log_{e}(1+\lambda)}{\lambda}}\right) = \left(\sqrt{\frac{\log_{e}(1+\lambda)}{\lambda}}\right)^{-n} H_{n}(x|\lambda).
\end{equation}
Further, it is to be noted that the $BVDHP$ $H_{n}(x,y|\lambda)$ has following representations in terms of degenerate Hermite polynomials $H_{n}(x|\lambda)$:
\begin{equation}\label{dheq9}
	H_{n}(x,y|\lambda) = (-i)^{n}y^{n/2}H_{n}\left(\dfrac{ix}{2\sqrt{y}} \bigg| \lambda\right)= (-y)^{n/2}H_{n}\left(\dfrac{x}{2\sqrt{-y}} \bigg| \lambda\right)= y^{n/2}H_{n}\left(\dfrac{x}{\sqrt{y}},1 \bigg| \lambda\right).
\end{equation}
This article explores new results that significantly advance the study of the $BVDHP$ $H_{n}(x,y|\lambda)$. By employing umbral approach, previously unexplored properties of these polynomials are established, leading to deeper insights and original contributions in this field. Our findings not only enhance the existing framework but also pave the way for further research through innovative methodologies.\\

The remainder of this paper is organized as follows. In Section $2$, the operational approach to the $BVDHP$ $H_{n}(x,y|\lambda)$ is considered. In Section $3$, certain results are derived using umbral methods. The monomiality and partial orthogonality properties of $BVDHP$ $H_{n}(x,y|\lambda)$ are established in Section $4$. In Section $5$, certain concluding remarks are considered.
\section{Operational approach}
The reliability and applicability of the proposed formalism are grounded in the theory of Hermite polynomials. These polynomials are well-known as solutions to the heat-type equations and play a crucial role in addressing problems arising in practical applications, see, for example, \cite{BDLS}.\\

Before proceeding to our main results, certain differential relations for the $BVDHP$ $H_{n}(x,y|\lambda)$ are derived in the following result:
\begin{prop}
	For the bivariate degenerate Hermite polynomials $H_{n}(x,y|\lambda)$, the following differential relations hold:
	\begin{equation}\label{dheq18}
		\partial_{x}^{r}H_{n}(x,y|\lambda) = \dfrac{n!}{(n-r)!}\left(\dfrac{\log(1+\lambda)}{\lambda}\right)^{r}H_{n-r}(x,y|\lambda),
	\end{equation}
	\begin{equation}\label{dheq19}
		\partial_{y}^{r}H_{n}(x,y|\lambda) = \dfrac{n!}{(n-2r)!}\left(\dfrac{\log(1+\lambda)}{\lambda}\right)^{r}H_{n-2r}(x,y|\lambda).
	\end{equation}
	\begin{proof}
		Differentiating equation \eqref{dheq3} w.r.t. $x$, we have
		\begin{equation*}
			\log(1+\lambda)\dfrac{t}{\lambda}(1+\lambda)^{(xt+yt^{2})/\lambda} = \sum_{n=0}^{\infty}\partial_{x}H_{n}(x,y|\lambda)\frac{t^{n}}{n!},
		\end{equation*}
	which in view of equation \eqref{dheq3} becomes
		\begin{equation}\label{dheq20}
			\dfrac{\log(1+\lambda)}{\lambda}\sum_{n=0}^{\infty}H_{n}(x,y|\lambda)\frac{t^{n+1}}{n!}=\sum_{n=0}^{\infty}\partial_{x}H_{n}(x,y|\lambda)\frac{t^{n}}{n!}.
		\end{equation}
		Equating the coefficients of $t^{n}$ in both side of equation \eqref{dheq20}, it follows that
		\begin{equation}\label{dheq21}
			\partial_{x}H_{n}(x,y|\lambda)=n\dfrac{\log(1+\lambda)}{\lambda}H_{n-1}(x,y|\lambda),
		\end{equation}
		which proves that relation \eqref{dheq18} is valid for $r=1$.\\
		Assuming that result \eqref{dheq18} is true for $r=k$, that is
		\begin{equation*}
			\partial_{x}^{k}H_{n}(x,y|\lambda)=\dfrac{n!}{(n-k)!}\left(\dfrac{\log(1+\lambda)}{\lambda}\right)^{k}H_{n-k}(x,y|\lambda), \quad k\leq n \in \mathbb{N}.
		\end{equation*}
		Differentiating above equation w.r.t. $x$ and using equation \eqref{dheq21} in the right-hand side of the resultant equation, it follows that
		\begin{equation*}
			\partial_{x}^{k+1}H_{n}(x,y|\lambda)=\dfrac{n!}{(n-(k+1))!}\left(\dfrac{\log(1+\lambda)}{\lambda}\right)^{k+1}\partial_{x}H_{n-(k+1)}(x,y|\lambda),
		\end{equation*}
	which shows that relation \eqref{dheq18} holds for $r=k+1$. Consequently by the principle of mathematical induction, assertion \eqref{dheq18} is proved.\\
		
		Proceeding on the same line as above, assertion \eqref{dheq19} can be proved.
	\end{proof}
\end{prop}
In the next result, the Rodrigue's type formula for $H_{n}(x,y|\lambda)$ is established.
\begin{thm}
	For the bivariate degenerate Hermite polynomials $H_{n}(x,y|\lambda)$, the following Rodrigue's type formula with respect to the variable $x$ holds true:
	\begin{equation}\label{dheq25}
		H_{n}(x,y|\lambda) =\dfrac{(-1)^{n/2}}{2^{n}}(1+\lambda)^{-x^{2}/4y\lambda}\;\dfrac{\partial^{n}}{\partial_{x}^{n}}\left((1+\lambda)^{x^{2}/4y\lambda}\right).
	\end{equation}
	\begin{proof}
Let
	\begin{equation}\label{dheq10}
		F(x,y,t|\lambda)=(1+\lambda)^{(xt+yt^{2})/\lambda},
	\end{equation}
which in view of generating function \eqref{dheq3} and applying Taylor's expansion formula
	\begin{equation*}
		F(x,y,t|\lambda)=\sum_{n=0}^{\infty}\left[\dfrac{\partial^{n}F(x,y,t|\lambda)}{\partial_{t}^{n}}\right]_{t=0}\;\dfrac{t^{n}}{n!},
	\end{equation*}
takes the form
	\begin{equation}\label{dheq35}
		F(x,y,t|\lambda)=\sum_{n=0}^{\infty}H_{n}(x,y|\lambda)\dfrac{t^{n}}{n!}=\sum_{n=0}^{\infty}\left[\dfrac{\partial^{n}}{\partial_{t}^{n}}(1+\lambda)^{(xt+yt^{2})/\lambda}\right]_{t=0}\dfrac{t^{n}}{n!}.
	\end{equation}
	Since, in view of relation \eqref{dheq9}, we have
	\begin{equation*}
		\sum_{n=0}^{\infty}H_{n}(x,y|\lambda)\dfrac{t^{n}}{n!}=\sum_{n=0}^{\infty}(-i)^{n}y^{n/2}H_{n}\left(\dfrac{ix}{2\sqrt{y}}\bigg|\lambda\right)\dfrac{t^{n}}{n!}.
	\end{equation*}
Consequently, in view of equations \eqref{dheq8} and \eqref{dheq35}, it follows that
\begin{equation}\label{dheq40}
	\sum_{n=0}^{\infty}\left[\dfrac{\partial^{n}}{\partial_{t}^{n}}(1+\lambda)^{(xt+yt^{2})/\lambda}\right]_{t=0}\dfrac{t^{n}}{n!}=(-1)^{n}y^{n/2}\left[\dfrac{\partial^{n}}{\partial_{t}^{n}}(1+\lambda)^{(ixt/\sqrt{y}-t^{2})/\lambda}\right]_{t=0}.
\end{equation}
Using relation \eqref{dheq40} in equation \eqref{dheq35} and equating the coefficients of $t^{n}$ in the resultant equation, we find
	\begin{equation}\label{dheq42}
		H_{n}(x,y|\lambda)=\left[\dfrac{\partial^{n}}{\partial_{t}^{n}}(1+\lambda)^{(xt+yt^{2})/\lambda}\right]_{t=0}=(-1)^{n}y^{n/2}\left[\dfrac{\partial^{n}}{\partial_{t}^{n}}(1+\lambda)^{(ixt/\sqrt{y}-t^{2})/\lambda}\right]_{t=0}.
	\end{equation}
Expressing equation \eqref{dheq42} in the following form:
	\begin{equation*}
		H_{n}(x,y|\lambda)=(-1)^{n}y^{n/2}\left[\dfrac{\partial^{n}}{\partial_{t}^{n}}\exp\left\{(ixt/\sqrt{y}-t^{2})\dfrac{\log(1+\lambda)}{\lambda}\right\}\right]_{t=0},
	\end{equation*}
which on further simplification becomes
	\begin{equation}\label{dheq1}
		H_{n}(x,y|\lambda)=(-1)^{n}y^{n/2}(1+\lambda)^{-x^{2}/4y\lambda}\left[\dfrac{\partial^{n}}{\partial_{t}^{n}}\exp\left\{-\left(\dfrac{ix}{2\sqrt{y}}-t\right)^{2}\dfrac{\log(1+\lambda)}{\lambda}\right\}\right]_{t=0}.
	\end{equation}
	For an arbitrary constant $A$, we have
	\begin{equation*}
		\dfrac{\partial}{\partial_{t}}f(Ax-t)=-A\dfrac{\partial}{\partial_{x}}f(Ax-t),
	\end{equation*}
consequently
	\begin{equation}\label{dheq68}
		\dfrac{\partial^{n}}{{\partial_{t}}^{n}}f(Ax-t)=(-1)^{n}A^{n}\dfrac{\partial^{n}}{{\partial_{x}}^{n}}f(Ax-t).
	\end{equation}
	Keeping $y$ fixed in equation \eqref{dheq1} and using differential relation \eqref{dheq68}, it follows that
	\begin{equation*}
		H_{n}(x,y|\lambda)=\dfrac{(-1)^{n/2}}{2^{n}}(1+\lambda)^{-x^{2}/4y\lambda}\left[\dfrac{\partial^{n}}{\partial_{x}^{n}}\exp\left(-\left(\dfrac{ix}{2\sqrt{y}}-t\right)^{2}\dfrac{\log(1+\lambda)}{\lambda}\right)\right]_{t=0},
	\end{equation*}
	which proves assertion \eqref{dheq25}.
	\end{proof}
\end{thm}
In order to find the expansion of $x^{n}$ in terms of the $BVDHP$ $H_{n}(x,y|\lambda)$, following result is proved:
\begin{thm}
	For the bivariate degenerate Hermite polynomials $H_{n}(x,y|\lambda)$, the following inverse explicit representation holds:
	\begin{equation}\label{dheq28}
		x^{n} = \sum_{r=0}^{\lfloor n/2\rfloor}\dfrac{(-1)^{r}n!}{r!(n-2r)!}\left(\dfrac{\log(1+\lambda)}{\lambda}\right)^{r-n}y^{r}H_{n-2r}(x,y|\lambda).
	\end{equation}
	\begin{proof}
	Rewriting operational relation \eqref{dheq23} as:
		\begin{equation}\label{dheq27}
		\left(\frac{x\log(1+\lambda)}{\lambda}\right)^{n}=\exp{\left(-\dfrac{\lambda}{\log(1+\lambda)}y\partial_{x}^{2}\right)}H_{n}(x,y|\lambda).
	\end{equation}
		Expanding the exponential in the left-hand side of equation \eqref{dheq27}, we find
		\begin{equation*}
		\left(\frac{x\log(1+\lambda)}{\lambda}\right)^{n}=	 \sum_{r=0}^{\infty}\dfrac{(-1)^{r}}{r!}\left(\dfrac{\lambda}{\log(1+\lambda)}\right)^{r}y^{r}\partial_{x}^{2r}H_{n}(x,y|\lambda),
		\end{equation*}
	which on using differential relation \eqref{dheq18} yields assertion \eqref{dheq28}.
	\end{proof}
\end{thm}
The operational approach adopted in establishing the above results appears to be well-suited for addressing the $BVDHP$ $H_{n}(x,y|\lambda)$.\\

Integral representations play a fundamental role in mathematical analysis, providing a powerful tool to express functions in terms of integrals. The importance of integral representations provides motivation to establish integral representations for the $BVDHP$ $H_{n}(x,y|\lambda)$.
\begin{prop}
	For the bivariate degenerate Hermite polynomials \( H_{n}(x,y|\lambda) \), the following integral representations hold:
	\begin{equation}\label{dheq47}
		\int_{0}^{x} H_{n}(\zeta,y|\lambda) d\zeta = \sum_{r=0}^{n} \dfrac{(-1)^{r}}{r+1} \binom{n}{r} (x)^{r+1} \left(\dfrac{\log(1+\lambda)}{\lambda}\right)^{r} H_{n-r}(x,y|\lambda),
	\end{equation}
	\begin{equation}\label{dheq48}
		\int_{0}^{y} H_{n}(x,\eta|\lambda) d\eta = \sum_{r=0}^{\left\lfloor\frac{n}{2}\right\rfloor} \dfrac{(-1)^{r}}{(r+1)!} (y)^{r+1} \dfrac{n!}{(n-2r)!} \left(\dfrac{\log(1+\lambda)}{\lambda}\right)^{r} H_{n-2r}(x,y|\lambda).
	\end{equation}
	\begin{proof}
	We know that the integral of a function \( g \in \mathbb{C}^{\infty} \) can be represented as a series:
		\begin{equation}\label{dheq46}
			\int_{0}^{x} g(\zeta) d\zeta = \sum_{r=0}^{\infty} (-1)^{r} \dfrac{x^{r+1}}{(r+1)!} g^{(r)}(x), \quad \forall x \in \mathbb{R},
		\end{equation}
		where \( g^{(r)}(x) \) signifies the r-th order derivative of the integrand function.\\
		
		Utilizing equation \eqref{dheq46} for $g(\zeta)=H_{n}(\zeta,y|\lambda)$, we have
		\begin{equation*}
			\int_{0}^{x} H_{n}(\zeta,y|\lambda) d\zeta = \sum_{r=0}^{n} \dfrac{(-1)^{r} x^{r+1}}{(r+1)!} \partial^{r}_{x} H_{n}(x,y|\lambda),
		\end{equation*}
		which in view of differential relation \eqref{dheq18} yields assertion \eqref{dheq47}.\\
		
		Following the same lines and using differential relation \eqref{dheq19}, integral representation \eqref{dheq48} can be obtained.
	\end{proof}
\end{prop}
In the next result, integrals involving the product of $H_{n}(x,y|\lambda)$ and $cosine$ functions are derived.
\begin{thm}
	For the product of bivariate degenerate Hermite polynomials \( H_{n}(x,y|\lambda) \) with $\text{cosine}$ function, the following integral representations hold:
	\begin{equation}\label{dheq51}
		\int_{0}^{x} H_{n}(\zeta,y|\lambda) \cos(\zeta) d\zeta = \sum_{r=0}^{n} (-1)^{r} \dfrac{\cos\left(x + r \dfrac{\pi}{2}\right)}{r+1} \binom{n}{r} \left(\dfrac{\log(1+\lambda)}{\lambda}\right)^{r} H_{n-r}(x,y|\lambda),
	\end{equation}
	\begin{equation}\label{dheq52}
		\int_{0}^{y} H_{n}(x,\eta|\lambda) \cos(\eta) d\eta = \sum_{r=0}^{\left\lfloor\frac{n}{2}\right\rfloor} (-1)^{r} \dfrac{\cos\left(x + r \dfrac{\pi}{2}\right)}{(r+1)!} \dfrac{n!}{(n-2r)!} \left(\dfrac{\log(1+\lambda)}{\lambda}\right)^{r} H_{n-2r}(x,y|\lambda).
	\end{equation}
	\begin{proof}
		For any \( f \in \mathbb{C}^{\infty} \) such that \( x \in \mathbb{R} \), the following relation holds:
		\begin{equation}\label{dheq49}
		{_0{D_{x}^{-1}}} (f(x) g(x))=	\int_{0}^{x} f(\zeta) g(\zeta) d\zeta,
		\end{equation}
		where the negative derivative operators ${_\beta{D_{x}^{-1}}}$ and ${_\beta{D_{x}^{-n}}}$ are defined as:
		\begin{equation*}
			{_\beta{D_{x}^{-1}}} h(x) = \int_{\beta}^{x} h(\zeta) d\zeta
		\end{equation*}
		and
		\begin{equation*}
		{_\beta{D_{x}^{-n}}}=	\underbrace{\int_{\beta}^{x} \cdots \int_{\beta}^{x}}_{n \text{ times}} h(\zeta) \, (d\zeta)^n =  \int_{\beta}^{x} \dfrac{(x-\zeta)^{n-1}}{(n-1)!}h(\zeta) \, d\zeta
		\end{equation*}
		respectively.\\
		By employing the generalized Leibniz formula
		\begin{equation}\label{dheq50}
			{_0{D_{x}^{-1}}} (f(x) g(x))=\sum_{r=0}^{\infty}\binom{-1}{r}f^{(r)}(x)\;{_0{D_{x}^{-r-1}}}[g(x)]=\sum_{r=0}^{\infty} (-1)^{r} f^{(r)}(x) g^{(-r-1)}(x)
		\end{equation}
		Using formula \eqref{dheq50} with $f(x)=H_{n}(x,y|\lambda)$ and $g(x)=\cos x$ in equation \eqref{dheq49} and simplifying, assertion \eqref{dheq51} is proved.\\
		
		Similarly, by keeping $y$ fixed and proceeding on the same lines as above with $f(y)=H_{n}(x,y|\lambda)$ and $g(y)=\cos y$, integral \eqref{dheq52} can be obtained.
	\end{proof}
\end{thm}

 The umbral formalism proves to be particularly advantageous for deriving novel properties of given special polynomials in a straightforward manner. In order to bolster the contention of using umbral methods for the $BVDHP$ $H_{n}(x,y|\lambda)$, these polynomials are reformulated using the umbral approach. In the next section, the umbral form of these polynomials is introduced and certain results are derived using the umbral methods.
 \section{Umbral approach}
 The foundation of umbral formalism lies in the classical Roman-Rota umbral calculus \cite{RT}, which establishes abstract notations and operational rules for manipulating formal power series. The essential formal tools are constructed using the ``image" and ``vacuum" functions. The broad spectrum of applications of umbral calculus continue to influence diverse area of mathematics, making it a subject of sustained research and interest, see, for example, \cite{DSML,LD,NU}.\\
 
 To establish the umbral representation of the $BVDHP$ $H_{n}(x,y|\lambda)$, we introduce the umbral operator $_\nu{\hat{h}_{y,\lambda}}$, defined as follows:
 \begin{equation*}
 	_\nu{\hat{h}_{y,\lambda}}^{r}\;:\;\phi(\nu)\to\phi(\nu+r).
 \end{equation*}
Here, the function $\phi(\nu)$ defined by
 \begin{equation}\label{dheq7}
 	\phi(\nu):=\phi_{\nu}=\frac{y^{\lfloor \nu/2\rfloor}}{\lfloor \nu/2\rfloor!}\left(\frac{\log(1+\lambda)}{\lambda}\right)^{\lfloor \nu/2\rfloor}\nu!\left|\cos\frac{\nu\pi}{2}\right|
 \end{equation}
 is referred to as the umbral vacuum, with the initial state $\phi(0)=\phi_{0}=1$.\\
 This operator satisfies the following property:
 \begin{equation*}
 	_\nu{\hat{h}_{y,\lambda}}^{m}\;{_\nu{\hat{h}_{y,\lambda}}^{n}}={_\nu{\hat{h}_{y,\lambda}}^{m+n}}.
 \end{equation*} 
We now proceed to formulate $BVDHP$ $H_{n}(x, y|\lambda)$ within an umbral framework. Consider
\begin{equation}\label{dheq4}
	(1+\lambda)^{yt^{2}/\lambda} =	e^{_\nu{\hat{h}_{y,\lambda}}t}\;\phi_{\nu}|_{\nu=0}.
\end{equation}
By expanding the exponential on the right-hand side and simplifying, we obtain
\begin{equation}\label{dheq6}
	(1+\lambda)^{yt^{2}/\lambda}=e^{\hat{h}_{y,\lambda}t}\phi_{0}=\sum_{n=0}^{\infty}\frac{y^{\lfloor n/2\rfloor}}{\lfloor n/2\rfloor!}\left(\frac{\log(1+\lambda)}{\lambda}\right)^{\lfloor n/2\rfloor}n!\dfrac{t^{n}}{n!}\left|\cos\frac{n\pi}{2}\right|,
\end{equation}
where $\hat{h}_{y,\lambda}$ is an umbral operator acting on the initial state $\phi_{0}$ such that
\begin{equation}\label{dheq5}
	\hat{h}_{y,\lambda}^{n}\phi_{0} = \phi_{n}=\frac{y^{\lfloor n/2\rfloor}}{\lfloor n/2\rfloor!}\left(\frac{\log(1+\lambda)}{\lambda}\right)^{\lfloor n/2\rfloor}n!\left|\cos\frac{n\pi}{2}\right|=\begin{cases}
		\dfrac{y^{s}(2s)!}{s!}\left(\dfrac{\log(1+\lambda)}{\lambda}\right)^{s}, & \text{if $r=2s$}, \\
		0, & \text{if $r=2s+1$},
	\end{cases}.
\end{equation}
Using equation \eqref{dheq6}, generating function \eqref{dheq3} can be expressed in the umbral form as:
\begin{equation*}
	\exp{(\hat{h}_{y,\lambda}t)}\exp{\left(\dfrac{x\log(1+\lambda)}{\lambda}t\right)}\phi_{0} = \sum_{n=0}^{\infty}H_{n}(x,y|\lambda)\frac{t^{n}}{n!}.
\end{equation*}
This leads to the umbral generating function of $BVDHP$ $H_{n}(x,y|\lambda)$ given by
\begin{equation}\label{dheq11}
	\exp{\left(\left(\hat{h}_{y,\lambda}+\dfrac{x\log(1+\lambda)}{\lambda}\right)t\right)}\phi_{0} = \sum_{n=0}^{\infty}H_{n}(x,y|\lambda)\frac{t^{n}}{n!}.
\end{equation}
Expanding the exponential in the left-hand side of generating equation \eqref{dheq11} and equating the coefficients of $t^{n}$, we find the following umbral representation of $BVDHP$ $H_{n}(x,y|\lambda)$:
\begin{equation}\label{dheq12}
	H_{n}(x,y|\lambda) = \left(\hat{h}_{y,\lambda} + \dfrac{x\log(1+\lambda)}{\lambda}\right)^{n}\phi_{0}.
\end{equation}
For $\lambda \to 0$, equations \eqref{dheq7}, \eqref{dheq5}, and \eqref{dheq12} reduce to the umbral vacuum, umbral operator, and umbral representation, respectively, associated with the bivariate Hermite polynomials $H_{n}(x,y)$, as introduced by Dattoli {\em et al.} \cite{DG}.\\

This limiting case reveals the intrinsic umbral structure of the $BVDHP$ $H_{n}(x,y|\lambda)$ and facilitates a transparent analysis of its umbral properties.\\

Additionally, differentiating \eqref{dheq11} $n$ times w.r.t. $t$ recursively, we derive
\begin{equation}\label{dheq13}
	\frac{\partial^{n}}{\partial t^{n}}\left[\exp{\left(\left(\hat{h}_{y,\lambda} + \dfrac{x\log(1+\lambda)}{\lambda}\right)t\right)}\right]\bigg|_{t=0}\phi_{0} = H_{n}(x,y|\lambda).
\end{equation}
Further, replacing $x$ by $2x$ and $y$ by $-1$ in umbral definition \eqref{dheq12}, we observe that
\begin{equation}\label{dheq16}
	H_{n}(2x, -1|\lambda)=H_{n}(x|\lambda)=\left(\dfrac{x\log(1+\lambda)}{\lambda}+\hat{h}_{\lambda}\right)\phi_{0},
\end{equation}
where
\begin{equation}\label{dheq22}
	\hat{h}_{\lambda}^{n}\phi_{0} = \frac{(-1)^{\lfloor n/2\rfloor}}{\lfloor n/2\rfloor!}\left(\frac{\log(1+\lambda)}{\lambda}\right)^{\lfloor n/2\rfloor}n!\left|\cos\frac{n\pi}{2}\right|.
\end{equation}

Expanding right-hand side of equation \eqref{dheq16} and using operator equation \eqref{dheq22} in the resultant, the series representation of degenerate Hermite polynomials (DHP) $H_{n}(x|\lambda)$ is obtained as:
\begin{equation}\label{dheq17}
	H_{n}(x|\lambda)=\sum_{k=0}^{\infty}\dfrac{(-1)^{k}n!}{k!(n-2k)!}\left(\dfrac{\log(1+\lambda)}{\lambda}\right)^{n-k}(2x)^{n-2k}.
\end{equation}
To gain a deeper understanding and validate the framework of the umbral formalism, we can easily verify series representation \eqref{dheq14} and heat equation \eqref{dheq15} of the $BVDHP$ $H_{n}(x,y|\lambda)$ by using umbral formalism, which aligns with results established in \cite{HR}.\\

Moreover, certain results involving Stirling numbers \cite{Young} and summation formulae, as mentioned in \cite{HR}, can be verified using the umbral formalism.\\

By employing standard calculus techniques along with the aforementioned formalism, we derive an infinite integral representation in the following Lemma:
\begin{lem}
	The following infinite integral involving $(1+\lambda)^{-yx^{2}/\lambda}$ holds:
	\begin{equation}\label{dheq54}
		\left|\int_{0}^{\infty} x^{\mu-1} (1+\lambda)^{-yx^{2}/\lambda} dx \right| = \dfrac{1}{2} \dfrac{1}{y^{\lfloor\mu/2\rfloor}} \left(\dfrac{\lambda }{\log(1+\lambda)}\right)^{\left\lfloor \dfrac{\mu}{2} \right\rfloor} \Gamma\left(\dfrac{\mu}{2}\right).
	\end{equation}
	\begin{proof}
		Consider the following identity:
		\begin{equation}\label{dheq53}
			e^{i\hat{h}_{y,\lambda}x}\phi_{0} = (1+\lambda)^{-yx^{2}/\lambda}.
		\end{equation}
		Using identity \eqref{dheq53}, the left-hand side of equation \eqref{dheq54} is expressed in the following form:
		\begin{equation*}
			\int_{0}^{\infty} x^{\mu-1} (1+\lambda)^{-yx^{2}/\lambda} dx = \int_{0}^{\infty} x^{\mu-1} e^{i\hat{h}_{y,\lambda}x} \phi_{0} dx.
		\end{equation*}
		Substituting $i\hat{h}_{y,\lambda}x = -t$ and utilizing the integral representation of the Gamma function, it follows that
		\begin{equation*}
			\int_{0}^{\infty} x^{\mu-1} (1+\lambda)^{-yx^{2}/\lambda} dx = \left(\dfrac{-1}{i\hat{h}_{y,\lambda}}\right)^{\mu}\phi_{0} \Gamma(\mu).
		\end{equation*}  
		Applying operator action \eqref{dheq5} in the right-hand side of the above equation, it simplifies to
		\begin{equation*}
			\int_{0}^{\infty} x^{\mu-1} (1+\lambda)^{-yx^{2}/\lambda} dx = \dfrac{(i)^{\mu} y^{\lfloor-\mu/2\rfloor} \Gamma(\mu) \Gamma(1-\mu)}{\Gamma\left(1-\dfrac{\mu}{2}\right)} \left(\dfrac{\log(1+\lambda)}{\lambda}\right)^{\lfloor-\mu/2\rfloor} \left| \cos \dfrac{\pi\mu}{2} \right|.
		\end{equation*}
		Utilizing the well-established reflection formula
		\begin{equation*}
			\Gamma(\nu) \Gamma(1-\nu) = \dfrac{\pi}{\sin(\pi\nu)},
		\end{equation*}
		and simplifying the resulting expression, we arrive at
		\begin{equation*}
			\int_{0}^{\infty} x^{\mu-1} (1+\lambda)^{-yx^{2}/\lambda} dx = \dfrac{(i)^{\mu}}{2} \dfrac{1}{y^{\lfloor\mu/2\rfloor}} \left(\dfrac{\lambda}{\log(1+\lambda)}\right)^{\lfloor\mu/2\rfloor} \Gamma\left(\dfrac{\mu}{2}\right).
		\end{equation*}
		Taking the modulus on both sides result \eqref{dheq54} is obtained.
	\end{proof}
\end{lem}
Now, let us consider the generating function involving even-indexed $BVDHP$, formulated using umbral representation \eqref{dheq12} as follows:
\begin{equation}\label{dheq55}
	\sum_{n=0}^{\infty}H_{2n}(x,y|\lambda)\dfrac{t^{n}}{n!}=\exp\left\{t\left(\hat{h}_{y,\lambda} + \dfrac{x\log(1+\lambda)}{\lambda}\right)^{2}\right\}\phi_{0}.
\end{equation}

By utilizing equation \eqref{dheq55}, another significant aspect of this framework emerges through the Gaussian integral identity \cite{BDLS}, as demonstrated in the following theorem:

\begin{thm}
	For the even-indexed bivariate degenerate Hermite polynomials $H_{2n}(x,y|\lambda)$, the following generating equation holds:
	\begin{equation}\label{dheq56}
		\sum_{n=0}^{\infty}H_{2n}(x,y|\lambda)\dfrac{t^{n}}{n!}=\dfrac{1}{\sqrt{1-(1+\lambda)^{4ty/\lambda}}}\exp\left\{\dfrac{t(1+\lambda)^{2x/\lambda}}{1-(1+\lambda)^{4ty/\lambda}}\right\},\quad \left|(1+\lambda)^{4ty/\lambda}\right|<1.
	\end{equation}
	\begin{proof}
		Applying the Gaussian integral identity, we have
		\begin{equation*}
			\exp\left(t\left(\hat{h}_{y,\lambda} + \dfrac{x\log(1+\lambda)}{\lambda}\right)^{2}\right)\phi_{0}=\dfrac{1}{\sqrt{\pi}}\int_{-\infty}^{\infty}\exp\left(-\zeta^{2}+2\sqrt{(t)}\left(\hat{h}_{y,\lambda} + \dfrac{x\log(1+\lambda)}{\lambda}\right)\zeta\right)\phi_{0}d\zeta.
		\end{equation*}
	Rewriting the above expression as:
		\begin{align}\label{dheq57}
			\notag &\hspace{-0.8cm}\exp\left(t\left(\hat{h}_{y,\lambda} + \dfrac{x\log(1+\lambda)}{\lambda}\right)^{2}\right)\phi_{0}=\dfrac{1}{\sqrt{\pi}}\int_{-\infty}^{\infty}\exp\left\{-\zeta^{2}\right\}\exp\left\{2\sqrt{(t)}\hat{h}_{y,\lambda}\zeta\right\}\\
			&\hspace{8cm} \exp\left\{2\sqrt{(t)} \dfrac{x\log(1+\lambda)}{\lambda}\zeta\right\}\phi_{0}d\zeta,
		\end{align}
	becomes
		\begin{equation*}
			\exp\{2\sqrt{(t)}\hat{h}_{y,\lambda}\zeta\}\phi_{0}=(1+\lambda)^{4t\zeta^{2}y/\lambda},
		\end{equation*}
		in equation \eqref{dheq57} yields
		\begin{align*}
			&\hspace{-1.8cm}\exp\left\{t\left(\hat{h}_{y,\lambda} + \dfrac{x\log(1+\lambda)}{\lambda}\right)^{2}\right\}\phi_{0}=\dfrac{1}{\sqrt{\pi}}\int_{-\infty}^{\infty}\exp\left\{-\zeta^{2}\right\}(1+\lambda)^{4t\zeta^{2}y/\lambda}\\
			&\hspace{7cm} \exp\left\{2\sqrt{(t)} \dfrac{x\log(1+\lambda)}{\lambda}\zeta\right\}d\zeta,
		\end{align*}
		which simplifies to
		\begin{align*}
			&\hspace{-1.4cm}\exp\left\{t\left(\hat{h}_{y,\lambda} + \dfrac{x\log(1+\lambda)}{\lambda}\right)^{2}\right\}\phi_{0} = \dfrac{1}{\sqrt{\pi}} \int_{-\infty}^{\infty} \exp\Bigg\{-\zeta^{2} \left(1 - \dfrac{4ty}{\lambda} \log(1+\lambda)\right) \\
			&\hspace{9.2cm} + 2\sqrt{t} \dfrac{x\log(1+\lambda)}{\lambda} \zeta \Bigg\} d\zeta.
		\end{align*}
		
		Evaluating the integral on the right-hand side using the formula
		\begin{equation}\label{dheq58}
			\int_{-\infty}^{\infty}\exp\{-ax^{2}+bx\}dx=\sqrt{\dfrac{\pi}{a}}\exp\left\{\dfrac{b^{2}}{4a}\right\},
		\end{equation}
		we have
		\begin{equation*}
			\exp\left\{t\left(\hat{h}_{y,\lambda} + \dfrac{x\log(1+\lambda)}{\lambda}\right)^{2}\right\}\phi_{0} = \dfrac{1}{\sqrt{1-(1+\lambda)^{4ty/\lambda}}}\exp\left\{\dfrac{t\dfrac{x^{2}}{\lambda^{2}}\left(\log(1+\lambda)\right)^{2}}{1-\dfrac{4ty}{\lambda}\log(1+\lambda)}\right\},
		\end{equation*}
		which on simplifying and using equation \eqref{dheq55} proves result \eqref{dheq56}.
	\end{proof}
\end{thm}
In the forthcoming section, the monomiality and partial orthogonality of the $BVDHP$ $H_{n}(x,y|\lambda)$ are considered.
\section{Monomiality and partial orthogonality}
Dattoli \cite{Dattoli} reformulated the concept introduced by Steffensen \cite{ST} and developed a novel theory termed as monomiality principle. According to this perspective, a family of polynomials \(\{g_{n}(x)\}_{n\in\mathbb{N}}\) is referred to as quasi-monomial, provided there exist a multiplicative operator \(\hat{\mathcal{M}}\) and a derivative operator \(\hat{\mathcal{P}}\), which satisfy the following relations:
\begin{equation}\label{dheq29}
	\hat{\mathcal{M}}\{g_{n}(x)\} = g_{n+1}(x),
\end{equation}
\begin{equation}\label{dheq30}
	\hat{\mathcal{P}}\{g_{n}(x)\} = ng_{n-1}(x),
\end{equation}
along with the commutation relation:
\begin{equation*}
	[\hat{\mathcal{P}}, \hat{\mathcal{M}}] = \hat{\mathcal{P}}\hat{\mathcal{M}} - \hat{\mathcal{M}}\hat{\mathcal{P}} = 1.
\end{equation*}

Certain additional properties follow, such as:
\begin{equation}\label{dheq31}
	\hat{\mathcal{M}}\hat{\mathcal{P}}\{g_{n}(x)\} = n g_{n}(x),
\end{equation}
which represents the differential equation satisfied by \(g_{n}(x)\).

Moreover, \(g_{n}(x)\) can be explicitly constructed as:
\begin{equation}\label{dheq32}
	g_{n}(x) = \hat{\mathcal{M}}^{n}\{1\}; \quad g_{0}(x) = 1.
\end{equation}
The monomiality of $BVDHP$ $H_{n}(x,y|\lambda)$ is established by proving the following result:
\begin{thm}
	 bivariate degenerate Hermite polynomials $H_{n}(x,y|\lambda)$ are quasi-monomial w.r.t. the following multiplicative and derivative operators:
	\begin{equation}\label{dheq33}
		\hat{\mathcal{M}}=\dfrac{\log(1+\lambda)}{\lambda}x+2y\partial_{x},
	\end{equation}
	and
	\begin{equation}\label{dheq34}
		\hat{\mathcal{P}}=\dfrac{\lambda}{\log(1+\lambda)}\partial_{x},
	\end{equation}
	respectively.
	\begin{proof}
		Differentiating generating equation \eqref{dheq3}, w.r.t. $t$, we find
		\begin{equation*}
			\dfrac{\log(1+\lambda)}{\lambda}(x+2yt)(1+\lambda)^{(xt+yt^{2})/\lambda} = \sum_{n=1}^{\infty}H_{n}(x,y|\lambda)\dfrac{nt^{n-1}}{n!},
		\end{equation*}
		which on again using generating function \eqref{dheq3} in the l.h.s. and simplifying becomes
		\begin{equation*}
			\dfrac{\log(1+\lambda)}{\lambda}x\sum_{n=0}^{\infty}H_{n}(x,y|\lambda)\dfrac{t^{n}}{n!}+2\dfrac{\log(1+\lambda)}{\lambda}y\sum_{n=1}^{\infty}nH_{n-1}(x,y|\lambda)\dfrac{t^{n}}{n!}=\sum_{n=0}^{\infty}H_{n+1}(x,y|\lambda)\dfrac{t^{n}}{n!}.
		\end{equation*}
		Equating the coefficient of $t^{n}$ in both side of above equation, the following recurrence relation is obtained:
		\begin{equation}\label{dheq36}
			\dfrac{\log(1+\lambda)}{\lambda}xH_{n}(x,y|\lambda)+2n\dfrac{\log(1+\lambda)}{\lambda}yH_{n-1}(x,y|\lambda)=H_{n+1}(x,y|\lambda), \quad n\ge1.
		\end{equation}
		Further, using equation \eqref{dheq21} in equation \eqref{dheq36}, it follows that
		\begin{equation*}
			\left(	\dfrac{\log(1+\lambda)}{\lambda}x+2y\partial_{x}\right)H_{n}(x,y|\lambda)=H_{n+1}(x,y|\lambda),
		\end{equation*}
		which in view of monomiality equation \eqref{dheq29} yields expression for multiplicative operator $\hat{\mathcal{M}}$.\\
		
		Assertion \eqref{dheq34} follows as a direct consequence of equations \eqref{dheq21} and \eqref{dheq30}.
	\end{proof}
\end{thm}
\begin{rem}
	Substituting expressions \eqref{dheq33} and \eqref{dheq34} of the multiplicative and derivative operators in monomiality equation \eqref{dheq31}, we deduce the following result:
\end{rem}
\begin{cor}
	The bivariate degenerate Hermite polynomials $H_{n}(x,y|\lambda)$ satisfy the following differential equation:
	\begin{equation}\label{dheq37}
		\left(x\partial_{x}+\dfrac{2\lambda}{\log(1+\lambda)}y\partial_{x}^{2}-n\right)H_{n}(x,y|\lambda)=0.
	\end{equation}
\end{cor}
\begin{rem}
	From  equations \eqref{dheq32} and \eqref{dheq33}, it follows that
	\begin{equation}\label{dheq38}
		\left(\dfrac{\log(1+\lambda)}{\lambda}x+2y\partial_{x}\right)^{n}\{1\}=H_{n}(x,y|\lambda),
	\end{equation}
	which verifies series expansion \eqref{dheq14} of $H_{n}(x,y|\lambda)$.
\end{rem}
\begin{rem}
	In view of equation \eqref{dheq38} and umbral form \eqref{dheq12}, we observe the following correspondence between umbral form and multiplicative operator of $BVDHP$ $H_{n}(x,y|\lambda)$:
	\begin{equation}\label{dheq39}
		\hat{\mathcal{M}} \leftrightarrow{} \left(\hat{h}_{y,\lambda} + \dfrac{x\log(1+\lambda)}{\lambda}\right).
	\end{equation}
\end{rem}
Next, we establish the orthogonality property of the degenerate Hermite polynomials $H_{n}(x|\lambda)$.\\
Since, the classical Hermite polynomials $H_{n}(x)$ are known to be orthogonal w.r.t. the weight function $w(x) = e^{-x^{2}}$ over the interval $(-\infty, \infty)$, with the following orthogonality conditions:
\begin{equation}\label{dheq60}
	\int_{-\infty}^{+\infty} e^{-x^{2}} H_{n}(x) H_{m}(x) \, dx = 0, \quad n \neq m,
\end{equation}
and
\begin{equation}\label{dheq61}
	\int_{-\infty}^{+\infty} e^{-x^{2}} \left[H_{n}(x)\right]^{2} \, dx = \sqrt{\pi} \, 2^{n} n!, \quad n = m,
\end{equation}
 The orthogonality property of Hermite polynomials makes it invaluable in quantum mechanics, probability theory, and other fields of mathematical physics \cite{ARF}.\\

The orthogonality property of the degenerate Hermite polynomials $H_{n}(x|\lambda)$ is established in the following result:
\begin{prop}
	The degenerate Hermite polynomials $H_{n}(x|\lambda)$ are orthogonal w.r.t. the weight function $w(x) = (1+\lambda)^{-x^{2}/\lambda}$ over $(-\infty, \infty)$. Specifically,
	\begin{equation}\label{dheq62}
		\int_{-\infty}^{+\infty} (1+\lambda)^{-x^{2}/\lambda} H_{n}(x|\lambda) H_{m}(x|\lambda) \, dx = 0, \quad if \quad m\not=n
	\end{equation}	
	and	
	\begin{equation}\label{dheq63}
		\int_{-\infty}^{\infty} (1+\lambda)^{-x^{2}/\lambda} \left[H_{n}(x|\lambda)\right]^{2} \, dx = \sqrt{\frac{\pi\lambda}{\log(1+\lambda)}} \, \log(1+\lambda)^{2^{n}n!/\lambda},\quad if \quad m=n.
	\end{equation}
	\begin{proof}
		Using equations \eqref{dheq59} in orthogonality relation \eqref{dheq60}, it follows that
		\begin{equation*}
			\sqrt{\frac{\log_{e}(1+\lambda)}{\lambda}} \int_{-\infty}^{+\infty} \exp\left(-x^{2} \frac{\log_{e}(1+\lambda)}{\lambda}\right) H_{n}\left(x\sqrt{\frac{\log_{e}(1+\lambda)}{\lambda}}\right) H_{m}\left(x\sqrt{\frac{\log_{e}(1+\lambda)}{\lambda}}\right) \, dx = 0.
		\end{equation*}			
		This simplifies to
		\begin{equation*}
			\int_{-\infty}^{+\infty} (1+\lambda)^{-x^{2}/\lambda} \left(\sqrt{\frac{\log_{e}(1+\lambda)}{\lambda}}\right)^{-m-n+1} H_{n}(x|\lambda) H_{m}(x|\lambda) \, dx = 0,
		\end{equation*}	
		which proves assertion \eqref{dheq62}.\\
		
		Similarly use of relation \eqref{dheq59} in \eqref{dheq61}, gives
		
		\begin{equation*}
			\int_{-\infty}^{+\infty} (1+\lambda)^{-x^{2}/\lambda} \left[H_{n}(x|\lambda)\right]^{2} \, dx = \left(\sqrt{\frac{\log(1+\lambda)}{\lambda}}\right)^{2n-1} \sqrt{\pi} \, 2^{n} n!,
		\end{equation*}
		
		which on simplification yields assertion \eqref{dheq63}.
	\end{proof}
\end{prop}
In the next result, we establish the partial orthogonality of the $BVDHP$ $H_{n}(x, y|\lambda)$ in the variable $x$, while the variable $y$ is held fixed.
\begin{thm}
	The bivariate degenerate Hermite polynomials $H_{n}(x,y|\lambda)$ are partially orthogonal in the variable $x$ w.r.t. the weight function $(1+\lambda)^{-x^{2}/4y\lambda}$, that is
	\begin{equation}\label{dheq64}
		\int_{-\infty}^{+\infty} (1+\lambda)^{-x^{2}/4y\lambda} H_{n}(x, y|\lambda) H_{m}(x, y|\lambda) \, dx = 0,\;y \in \mathbb{C}\;\;if\;\;m\not=n
	\end{equation}
	and
	\begin{equation}\label{dheq65}
		\int_{-\infty}^{+\infty} (1+\lambda)^{-x^{2}/4y\lambda} \left[H_{n}(x, y|\lambda)\right]^{2} \, dx = \sqrt{\frac{\pi\lambda}{\log(1+\lambda)}} (-y)^{n+1/2}\, \log(1+\lambda)^{2^{n+1}n!/\lambda},\;y \in \mathbb{C}\;\;if\;\;m=n.
	\end{equation}
	\begin{proof}
		Since
		\begin{equation*}
			(1+\lambda)^{(2axt - a^{2}t^{2})/\lambda} = \sum_{n=0}^{\infty} H_{n}(x|\lambda) \frac{(at)^{n}}{n!},
		\end{equation*}
		it follows that
		\begin{equation*}
			H_{n}(x|\lambda) = a^{-n} H_{n}(2ax, -a^{2}|\lambda).
		\end{equation*}
		Since $H_{n}(2ax, -a^{2}|\lambda)$ is also orthogonal w.r.t. weight function $(1+\lambda)^{-x^{2}/\lambda}$. Therefore, we have
		\begin{equation}\label{dheq66}
			\int_{-\infty}^{+\infty} (1+\lambda)^{-x^{2}/\lambda} H_{n}(2ax, -a^{2}|\lambda) H_{m}(2ax, -a^{2}|\lambda) \, dx = 0,
		\end{equation}
		and
		\begin{equation}\label{dheq67}
			\int_{-\infty}^{+\infty} (1+\lambda)^{-x^{2}/\lambda} \left[H_{n}(2ax, -a^{2}|\lambda)\right]^{2} \, dx = \sqrt{\pi} \, 2^{n} a^{2n} n!.
		\end{equation}
		Substituting $2ax = u$ and $a^{2} = v$ in relation \eqref{dheq66} proves assertion \eqref{dheq64}. Similarly, using the same substitutions in relation \eqref{dheq67} proves \eqref{dheq65}.
	\end{proof}
\end{thm}
The problem associated with the possibility of defining negative order $BVDHP$ within the framework of umbral formalism presents a significant and promising research direction. In the next section, an idea is provided to explore further insights of this problem.
\section{Concluding remarks}
In this article, the previously unexplored properties of $BVDHP$ $H_{n}(x,y|\lambda)$ are investigated. The umbral form of these polynomials is introduced and generating equations and operational identities for these polynomials are established. The study also examines the monomiality principle and orthogonality property of the $BVDHP$ $H_{n}(x,y|\lambda)$.\\

To further extend the concept of umbral formalism, let us consider the Gaussian integral:
\begin{equation*}
	I(a,b|\lambda) = \int_{-\infty}^{\infty} \exp\left\{-x^{2}\left(a\dfrac{\log(1+\lambda)}{\lambda} + \hat{h}_{b,\lambda}\right)\right\} \phi_{0} \;dx,
\end{equation*}
which on using integral formula \eqref{dheq58} gives
\begin{equation}\label{dheq45}
	I(a,b|\lambda) = \sqrt{\dfrac{\pi}{\hat{H}(a,b|\lambda)}} \phi_{0} = \sqrt{\pi} \hat{H}_{-1/2}(a,b|\lambda) \phi,
\end{equation}
where
\begin{equation*}
	\hat{H}(a,b|\lambda) = \left(a \dfrac{\log(1+\lambda)}{\lambda} + \hat{h}_{b,\lambda}\right).
\end{equation*}

This result hints at the possibility of defining negative-order Hermite polynomials. Here, we introduce the negative-order degenerate Hermite functions (NODHF) $H_{-\mu}(x,y|\lambda)$, inspired by the umbral form \eqref{dheq12} as follows:
\begin{equation}\label{dheq43}
	H_{-\mu}(x,y|\lambda) = \left(\hat{h}_{y,\lambda} + \dfrac{x\log(1+\lambda)}{\lambda}\right)^{-\mu} \phi_{0}, \quad \forall \mu \in \mathbb{R}^{+}, \forall x, y \in \mathbb{R}.
\end{equation}

These are no longer polynomials but instead represent Hermite functions.\\

 Further, employing the Laplace transform, it follows that
\begin{equation*}
	H_{-\mu}(x,y|\lambda) = \dfrac{1}{\left(\hat{h}_{y,\lambda} + \dfrac{x\log(1+\lambda)}{\lambda}\right)^{\mu}} \phi_{0} = \dfrac{1}{\Gamma(\mu)} \int_{0}^{\infty} \zeta^{\mu-1} e^{-x\zeta} e^{-\hat{h}_{x,\lambda}\zeta} \phi_{0} d\zeta,
\end{equation*}
which on utilizing equation \eqref{dheq6}, yields
\begin{equation}\label{dheq44}
	H_{-\mu}(x,y|\lambda) = \dfrac{1}{\Gamma(\mu)} \int_{0}^{\infty} \zeta^{\mu-1} e^{-x\zeta} e^{-y\zeta^{2}/\lambda\;(\log(1+\lambda))} d\zeta.
\end{equation}

The above framework opens new possibilities for addressing practical problems and may serve as a foundation for future research in this area.\\

{\bf Declaration of competing interest}\\

None.


\begin{thebibliography}{99}
	\bibitem{ARF} G. B. Arfken, H. J. Weber, F. E. Harris {\em Mathematical Methods for Physicists} Academic Press (2013).
	\bibitem{APP} P. Appell, J. Kamp\'e de F\'eriet, {\em Fonctions Hyperg\'eom\'etriques et Hypersph\'eriques: Polyn$\hat{o}$mes d’ Hermite}, Gauthier-Villars, Paris (1926).
\bibitem{BDLS}	D. Babusci, G. Dattoli, S. Licciardi, E. Sabia, {\em Mathematical Methods for Physicists}, Singapore, World Scientific Publishing (2019).
	\bibitem{Dattoli} G. Dattoli, {\em Generalized polynomials, operational identities and their applications}, J. Comput. Appl. Math. {\bf 118} 111--123 (2000).
	\bibitem{DG} G. Dattoli, B. Germano, M. R. Martinelli, P. E. Ricci, {\em Lacunary generating functions of Hermite polynomials and symbolic methods}, Ilirias J. Math. {\bf 4} 16–23 (2015).
	\bibitem{DSML} G. Dattoli, Subuhi Khan, Mehnaz Haneef, S. Licciardi, {\em On umbral property of a family of hyperbolic-like functions appearing in magnetic transport problem}, Rep. Math. Phys. {\bf 92}(1) 37--48 (2023).
\bibitem{HR} K. W. Hwang, C. S. Ryoo, {\em  Differential equations associated with two variable degenerate Hermite polynomials}, Mathematics {\bf 8} (2) 228 (2020).
\bibitem{KJLK} T. Kim, D. S. Kim, L. C. Jang, H. Lee, H. Kim, {\em Representations of degenerate Hermite polynomials}, Adv. Appl. Math. {\bf 139} (3) 102359 (2022).
\bibitem{LD}S. Licciardi, G. Dattoli, {\em Guide to the Umbral Calculus, A Different Mathematical Language}, World Scientific (2022).
\bibitem{NU} N. Raza, U. Zainab, S. Araci, A. Esi, {\em Identities involving 3-variable Hermite polynomials arising from umbral method}, Adv. Differ. Equ. 640, 1–16 (2020).
	\bibitem{RT}S. M. Roman, G. C. Rota, {\em The umbral calculus}, Adv. in Math. {\bf 27}, 95--188 (1978).
	\bibitem{ST}J. F. Steffensen, {\em The poweroid, an extension of the mathematical notion of power}, Acta Math. {\bf 73}, 333--366 (1941).
\bibitem{Young} P. T. Young, {\em Degenerate Bernoulli polynomials, generalized factorial sums, and their applications}, J. Number Theorey  {\bf128} 738–758 (2008).

\end{thebibliography}
\end{document}